\documentclass[12pt]{article}

\usepackage{graphicx}
\usepackage{amssymb}
\usepackage{amsfonts}
\usepackage{moreverb}

\begin{document}

\begin{center}
\Large\bf GloptiPoly 3: moments, optimization and semidefinite
programming
\end{center}

\begin{center}
Didier Henrion$^{1,2}$, Jean-Bernard Lasserre$^{1,3}$, Johan L\"ofberg$^4$
\end{center}

\footnotetext[1]{LAAS-CNRS, University of Toulouse, France}

\footnotetext[2]{Faculty of Electrical Engineering, Czech Technical University in
Prague, Czech Republic}

\footnotetext[3]{Institute of Mathematics, University of Toulouse, France}

\footnotetext[4]{Department of Electrical Engineering, Link\"oping University,
Sweden}

\begin{center}
Version 3.0 of \today
\end{center}

\begin{abstract}
We describe a major update of our Matlab freeware GloptiPoly
for parsing generalized problems of moments and solving them
numerically with semidefinite programming.
\end{abstract}

\section{What is GloptiPoly ?}

Gloptipoly 3 is intended to solve, or at least approximate,
the Generalized Problem of Moments (GPM), an infinite-dimensional
optimization problem which can be viewed as an extension of
the classical problem of moments \cite{gpm}.
From a theoretical viewpoint, the GPM has developments and impact
in various areas of mathematics such as algebra, Fourier analysis,
functional analysis, operator theory, probability and statistics,
to cite a few. In addition, and despite a rather simple and
short formulation, the GPM has a large number of important
applications in various fields such as optimization, probability,
finance, control, signal processing, chemistry, cristallography,
tomography, etc. For an account of various methodologies as well
as some of potential applications, the interested reader is referred to
\cite{akhiezer,akhiezer2} and the nice collection of papers \cite{landau}.

The present version of GloptiPoly 3 can handle moment problems with 
polynomial data. Many important applications in e.g. optimization,
probability, financial economics and optimal control, can be viewed
as particular instances of the GPM, and (possibly after some
transformation) of the GPM with polynomial data. 

The approach is similar to that used in the former version 2 of
GloptiPoly \cite{gloptipoly2}. The software allows to build up a
hierarchy of semidefinite programming (SDP), or linear matrix
inequality (LMI) relaxations of the GPM, whose associated monotone
sequence of optimal values converges to the global optimum.
For more details on the approach, the interested reader is
referred to \cite{gpm}.

\section{Installation}

GloptiPoly 3 is a freeware subject to the General Public Licence
(GPL) policy. It can be downloaded at 
\begin{center}
{\tt www.laas.fr/$\sim$henrion/software/gloptipoly3}
\end{center}

The package, available as a compressed archive, consists of several
m-files and subdirectories, and it contains no binaries. Extracted
files are placed in a {\tt gloptipoly3} directory that should be
declared in the Matlab working path, using e.g. Matlab's command
\begin{verbatim}
>> addpath gloptipoly3
\end{verbatim}

GloptiPoly 3 uses by default the semidefinite programming solver
SeDuMi \cite{sedumi}, so this package should be properly installed.
Other semidefinite solvers can also be used provided they are
installed and interfaced through YALMIP \cite{yalmip}.

\section{Getting started}\label{started}

Please type the command
\begin{verbatim}
>> gloptipolydemo
\end{verbatim}
to run interactively the basic example that follows.

Consider the classical problem of minimizing globally the two-dimensional
six-hump camel back function \cite{gloptipoly2}
\[
\min_{x \in {\mathbb R}^2} g_0(x) = 
4x_1^2+x_1x_2-4x_2^2-2.1x_1^4+4x_2^4+\frac{1}{3}x_1^6.
\]
The function has six local minima, two of them being global minima.

Using GloptiPoly 3, this optimization problem can be modeled as a moment
problem as follows: 
\begin{verbatim}
>> mpol x1 x2
>> g0 = 4*x1^2+x1*x2-4*x2^2-2.1*x1^4+4*x2^4+x1^6/3
Scalar polynomial
4x1^2+x1x2-4x2^2-2.1x1^4+4x2^4+0.33333x1^6
>> P = msdp(min(g0));
GloptiPoly 3.0
Define moment SDP problem
...
(GloptiPoly output suppressed)
...
Generate moment SDP problem
>> P
Moment SDP problem 
  Measure label             = 1
  Relaxation order          = 3
  Decision variables        = 27
  Semidefinite inequalities = 10x10
\end{verbatim}

Once the moment problem is modeled, a semidefinite solver can be used to
solve it numerically. Here we use SeDuMi
\cite{sedumi} which is assumed to be installed and accessible
from the Matlab working path:
\begin{verbatim}
>> [status,obj] = msol(P)                       
GloptiPoly 3.0
Solve moment SDP problem
*****************************************************
Calling SeDuMi
SeDuMi 1.1R3 by AdvOL, 2006 and Jos F. Sturm, 1998-2003.
...
(SeDuMi output suppressed)
...
2 globally optimal solutions extracted
>> status
status =
     1
>> obj
obj =
    -1.0316
>> x = double([x1 x2]);
x(:,:,1) =
    0.0898   -0.7127
x(:,:,2) =
   -0.0898    0.7127
\end{verbatim}
The flag {\tt status = 1} means that the moment problem is solved
successfully and that GloptiPoly can extract two globally optimal
solutions reaching the objective function {\tt obj = -1.0316}.

\section{From version 2 to version 3}

The major changes incorporated into GloptiPoly
when passing from version 2 to 3 can be summarized as follows:
\begin{itemize}
\item Use of native polynomial objects and object-oriented programming
with specific classes for multivariate polynomials, measures, moments,
and corresponding overloaded operators. In contrast with version 2,
the Symbolic Toolbox for Matlab (gateway to the Maple kernel)
is not required anymore to process polynomial data.
\item Generalized problems of moments featuring several measures with
semialgebraic support constraints and linear moment constraints
can be processed and solved. Version 2 was limited to moment problems
on a unique measure without moment constraints.
\item Explicit moment substitutions are carried out to reduce the
number of variables and constraints.
\item The moment problems can be solved numerically with any semidefinite
solver, provided it is interfaced through YALMIP. In contrast, version 2
used only the solver SeDuMi.
\end{itemize}

\section{Solving generalized problems of moments}\label{solving}

GloptiPoly 3 uses advanced Matlab features for object-oriented
programming and overloaded operators. The user should be familiar
with the following basic objects.

\subsection{Multivariate polynomials ({\tt mpol})}

A multivariate polynomial is an affine combination of monomials, each
monomial depending on a set of variables. Variables can be declared
in the Matlab working space as follows:
\begin{verbatim}
>> clear
>> mpol x
>> x
Scalar polynomial
x
>> mpol y 2
>> y
2-by-1 polynomial vector
(1,1):y(1)
(2,1):y(2)
>> mpol z 3 2
>> z
3-by-2 polynomial matrix
(1,1):z(1,1)
(2,1):z(2,1)
(3,1):z(3,1)
(1,2):z(1,2)
(2,2):z(2,2)
(3,2):z(3,2)
\end{verbatim}
Variables, monomials and polynomials are defined as objects of
class {\tt mpol}.

All standard Matlab operators have been overloaded for
{\tt mpol} objects:
\begin{verbatim}
>> y*y'-z'*z+x^3
2-by-2 polynomial matrix
(1,1):y(1)^2-z(1,1)^2-z(2,1)^2-z(3,1)^2+x^3
(2,1):y(1)y(2)-z(1,1)z(1,2)-z(2,1)z(2,2)-z(3,1)z(3,2)+x^3
(1,2):y(1)y(2)-z(1,1)z(1,2)-z(2,1)z(2,2)-z(3,1)z(3,2)+x^3
(2,2):y(2)^2-z(1,2)^2-z(2,2)^2-z(3,2)^2+x^3
\end{verbatim}

Use the instruction
\begin{verbatim}
>> mset clear
\end{verbatim}
to delete all existing GloptiPoly variables
from the Matlab working space.

\subsection{Measures ({\tt meas})}

Variables can be associated with real-valued measures, and one
variable is associated with only one measure.
For GloptiPoly, measures are identified with a label, a positive
integer. When starting a GloptiPoly session, the default measure
has label 1. By default, all created variables are associated
with the current measure. Measures can be handled with the
class {\tt meas} as follows:
\begin{verbatim}
>> mset clear
>> mpol x
>> mpol y 2 
>> meas
Measure 1 on 3 variables: x,y(1),y(2)
>> meas(y) % create new measure
Measure 2 on 2 variables: y(1),y(2)
>> m = meas 
1-by-2 vector of measures
1:Measure 1 on 1 variable: x
2:Measure 2 on 2 variables: y(1),y(2)
>> m(1)
Measure number 1 on 1 variable: x
\end{verbatim}
The above script creates a measure $d\mu_1(x)$ on $\mathbb R$ and a
measure $d\mu_2(y)$ on ${\mathbb R}^2$.

Use the instruction 
\begin{verbatim}
>> mset clearmeas
\end{verbatim}
to delete all existing GloptiPoly measures from
the working space. Note that this does not delete
existing GloptiPoly variables.

\subsection{Moments ({\tt mom})}\label{mom}

Linear combinations of moments of a given measure
can be manipulated with the {\tt mom} class as follows:
\begin{verbatim}
>> mom(1+2*x+3*x^2)
Scalar moment
I[1+2x+3x^2]d[1]
>> mom(y*y') 
2-by-2 moment matrix
(1,1):I[y(1)^2]d[2]
(2,1):I[y(1)y(2)]d[2]
(1,2):I[y(1)y(2)]d[2]
(2,2):I[y(2)^2]d[2]
\end{verbatim}
The notation {\tt I[p]d[k]} stands for $\int p \:d\mu_k$ where $p$ is a
polynomial of the variables associated with measure $d\mu_k$,
and $k$ is the measure label.

Note that it makes no sense to define moments over several measures,
or nonlinear moment expressions:
\begin{verbatim}
>> mom(x*y(1))
??? Error using ==> mom.mom
Invalid partitioning of measures in moments
>> mom(x)*mom(y(1))
??? Error using ==> mom.times
Invalid moment product
\end{verbatim}

Note also the distinction between a constant term and the mass of a
measure:
\begin{verbatim}
>> 1+mom(x)
Scalar moment
1+I[x]d[1]
>> mom(1+x)
Scalar moment
I[1+x]d[1]
>> mass(x)
Scalar moment
I[1]d[1]
\end{verbatim}

Finally, let us mention three equivalent notations
to refer to the mass of a measure:
\begin{verbatim}
>> mass(meas(y))
Scalar moment
I[1]d[2]
>> mass(y)
Scalar moment
I[1]d[2]
>> mass(2)
Scalar moment
I[1]d[2]
\end{verbatim}
The first command refers explicitly to the measure,
the second command is a handy short-cut to refer to a measure
via its variables, and the third command refers to GloptiPoly's
labeling of measures.

\subsection{Support constraints ({\tt supcon})}\label{supcon}

By default, a measure on $n$ variables is defined on the whole
${\mathbb R}^n$. We can restrict
the support of a mesure to a given semialgebraic set
as follows:
\begin{verbatim}
>> 2*x^2+x^3 == 2+x
Scalar measure support equality
2x^2+x^3 == 2+x
>> disk = (y'*y <= 1) 
Scalar measure support inequality
y(1)^2+y(2)^2 <= 1
\end{verbatim}
Support constraints are modeled by objects of class {\tt
supcon}. The first command means that variable $x$ must satisfy
$x^3+2x^2-x-2 = (x-1)(x+1)(x+2) = 0$, i.e. measure $d\mu_1(x)$ must
be discrete, a linear combination of three Dirac at
$1$, $-1$ and $-2$. The second command restricts measure
$d\mu_2(y)$ within the unit disk.

Note that it makes no sense to define a support constraint on
several measures:
\begin{verbatim}
>> x+y(1) <= 1
??? Error using ==> supcon.supcon
Invalid reference to several measures
\end{verbatim}

\subsection{Moment constraints ({\tt momcon})}

We can constrain linearly the moments of several
measures:
\begin{verbatim}
>> mom(x^2+2) == 1+mom(y(1)^3*y(2))
Scalar moment equality constraint
I[2+x^2]d[1] == 1+I[y(1)^3y(2)]d[2]
>> mass(x)+mass(y) <= 2
Scalar moment inequality constraint
I[1]d[1]+I[1]d[2] <= 2
\end{verbatim}
Moment constraints are modeled by objects of class {\tt
momcon}.

For GloptiPoly an objective function to be minimized or maximized
is considered as a particular moment constraint:
\begin{verbatim}
>> min(mom(x^2+2))
Scalar moment objective function
min I[2+x^2]d[1]
>> max(x^2+2)
Scalar moment objective function
max I[2+x^2]d[1]
\end{verbatim}
The latter syntax is a handy short-cut which directly converts
an {\tt mpol} object into an {\tt momcon} object.

\subsection{Floating point numbers ({\tt double})}\label{double}

Variables in a measure can be assigned numerical values:
\begin{verbatim}
>> m1 = assign(x,2)
Measure 1 on 1 variable: x
supported on 1 point
\end{verbatim}
which is equivalent to enforcing a discrete support for
the measure. Here $d\mu_1$ is set to the Dirac
at the point $2$.

The {\tt double} operator converts a measure or its
variables into a floating point number:
\begin{verbatim}
>> double(x)
ans =
     2
>> double(m1)
ans =
     2
\end{verbatim}

Polynomials can be evaluated similarly:
\begin{verbatim}
>>double(1-2*x+3*x^2)
ans =
     9
\end{verbatim}

Discrete measure supports consisting of several points
can be specified in an array:
\begin{verbatim}
>> m2 = assign(y,[-1 2 0;1/3 1/4 -2])
Measure 2 on 2 variables: y(1),y(2)
supported on 3 points
>> double(m2)
ans(:,:,1) =
   -1.0000
    0.3333
ans(:,:,2) =
    2.0000
    0.2500
ans(:,:,3) =
     0
    -2
\end{verbatim}

\subsection{Moment SDP problems ({\tt msdp})}

GloptiPoly 3 can manipulate and solve Generalized Problems of Moments (GPM)
as defined in \cite{gpm}:
\[
\begin{array}{ll}
\min_{d\mu} \:(\mathrm{or}\: \max) & \displaystyle
\sum_k \int_{{\mathbb K}_k} g_{0k}(x)
d\mu_k(x) \\
\mathrm{s.t.} & \displaystyle
\sum_k \int_{{\mathbb K}_k} h_{jk}(x) d\mu_k(x)
\geq\:(\mathrm{or}\: =)\: b_j, \quad j=0,1,\ldots \\
\end{array}
\]
where measures
$d\mu_k$ are supported on basic semialgebraic sets
\[
{\mathbb K}_k = \{x \in {\mathbb R}^{n_k} \: :\: g_{ik}(x) \geq 0,
\quad i=1,2\ldots\}.
\]
In the above notations, $g_{ik}(x)$, $h_{jk}(x)$ are given real polynomials
and $b_j$ are given real constants. The decision variables in the GPM
are measures $d\mu_k(x)$, and GloptiPoly 3 allows to optimize over them
through their moments
\[
y_{\alpha_k} = \int_{{\mathbb K}_k} x^{\alpha_k} d\mu_k(x), \quad
\alpha_k \in {\mathbb N}^{n_k}
\]
where the $\alpha_k$ are multi-indices.

\subsection{Solving moment problems {\tt msol}}\label{msol}

Once a moment problem is defined, it can be solved numerically
with the instruction {\tt msol}.
In the sequel we give several examples of GPMs
handled with GloptiPoly 3.

\subsubsection{Unconstrained minimization}\label{unconstrained}

Following \cite{lasserre}, given a multivariate polynomial $g_0(x)$,
the unconstrained optimization problem
\[
\min_{x \in {\mathbb R}^n} g_0(x)
\]
can be formulated as a linear moment optimization problem
\[
\begin{array}{ll}
\min_{d\mu} & \int g_0(x)d\mu(x) \\
\mathrm{s.t.} & \int d\mu(x) = 1
\end{array}
\]
where measure $d\mu$ lives in the space 
${\mathbb B}^n$ of finite Borel signed measures
on ${\mathbb R}^n$. The equality constraint
indicates that the mass of $d\mu$ is equal to one, or equivalently,
that $d\mu$ is a probability measure.

In general, this linear (hence convex) reformulation 
of a (typically nonconvex) polynomial problem is
not helpful because there is no computationally efficient way to
represent measures and their underlying Borel spaces.
The approach proposed in \cite{lasserre} consists in using convex
semidefinite representations of the space ${\mathbb B}^n$ truncated to
finite degree moments. GloptiPoly 3 allows
to input such moment optimization problems in an user-friendly way,
and to solve them using existing software for semidefinite
programming (SDP).

In Section \ref{started} we already encountered an example
of an unconstrained polynomial optimization solved with GloptiPoly 3.
Let us revisit this example:
\begin{verbatim}
>> mset clear
>> mpol x1 x2
>> g0 = 4*x1^2+x1*x2-4*x2^2-2.1*x1^4+4*x2^4+x1^6/3
Scalar polynomial
4x1^2+x1x2-4x2^2-2.1x1^4+4x2^4+0.33333x1^6
>> P = msdp(min(g0));
...
>> msol(P)
...
2 globally optimal solutions extracted
Global optimality certified numerically
\end{verbatim}
This indicates that the global minimum is attained
with a discrete measure supported on two points. The measure
can be constructed from the knowledge of its first moments
of degree up to 6:
\begin{verbatim}
>> meas
Measure 1 on 2 variables: x1,x2
  with moments of degree up to 6, supported on 2 points
>> double(meas)
ans(:,:,1) =
    0.0898
   -0.7127
ans(:,:,2) =
   -0.0898
    0.7127
>> double(g0)
ans(:,:,1) =
   -1.0316
ans(:,:,2) =
   -1.0316
\end{verbatim}

When converting to floating point numbers with the operator 
{\tt double}, it is essential to make the distinction between {\tt mpol}
and {\tt mom} objects:
\begin{verbatim}
>> v = mmon([x1 x2],2)'
1-by-6 polynomial vector
(1,1):1
(1,2):x1
(1,3):x2
(1,4):x1^2
(1,5):x1x2
(1,6):x2^2
>> double(v)
ans(:,:,1) =
    1.0000    0.0898   -0.7127    0.0081   -0.0640    0.5079
ans(:,:,2) =
    1.0000   -0.0898    0.7127    0.0081   -0.0640    0.5079
>> double(mom(v))
ans =
    1.0000    0.0000   -0.0000    0.0081   -0.0640    0.5079
\end{verbatim}
The first instruction {\tt mmon} generates a vector of monomials {\tt v}
of class {\tt mpol}, so the command {\tt double(v)} calls the convertor
{\tt @mpol/double} which evaluates a polynomial expression on
the discrete support of a measure (here two points).
The last command {\tt double(mom(v))} calls the convertor
{\tt @mom/double} which returns the value of the moments
obtained after solving the moment problem.

Note that when inputing moment problems on a unique measure
whose mass is not constrained,
GloptiPoly assumes by default that the measure has mass one,
i.e. that we are seeking a probability measure.
Therefore, if {\tt g0} is the polynomial defined previously,
the two instructions
\begin{verbatim}
>> P = msdp(min(g0));
\end{verbatim}
and
\begin{verbatim}
>> P = msdp(min(g0), mass(meas(g0))==1);
\end{verbatim}
are equivalent. See also Section \ref{mom} for handling
masses of measures and Section \ref{constrained} for
more information on mass constraints.

\subsubsection{Constrained minimization}\label{constrained}

Consider now the constrained polynomial optimization problem
\[
\min_{x \in {\mathbb K}} g_0(x)
\]
where
\[
{\mathbb K} = \{ x \in {\mathbb R}^n \: :\:
g_i(x) \geq 0, \:i=1,2,\ldots\}
\]
is a basic semialgebraic set described by given
polynomials $g_i(x)$.
Following \cite{lasserre}, this
(nonconvex polynomial) problem can be formulated as the
(convex linear) moment problem
\[
\begin{array}{ll}
\min_{d\mu} & \int_{\mathbb K} g_0(x)d\mu(x) \\
\mathrm{s.t.} & \int_{\mathbb K} d\mu(x) = 1
\end{array}
\]
where the indeterminate is a probability measure $d\mu$ of
${\mathbb B}^n$ which is
now supported on set $\mathbb K$. In other words
\[
\int_{{\mathbb R}^n/{\mathbb K}} d\mu(x) = 0.
\]

As an example, consider the non-convex quadratic problem
of Section 4.4 in \cite{gloptipoly2}:
\[
\begin{array}{ll}
\min & -2x_1+x_2-x_3 \\
\mathrm{s.t.} & 24-20x_1+9x_2-13x_3+4x_1^2-4x_1x_2+4x_1x_3+2x_2^2-2x_2x_3+2x_3^2 \geq 0 \\
& x_1+x_2+x_3 \leq 4, \quad 3x_2+x_3 \leq 6\\
& 0 \leq x_1 \leq 2, \quad 0 \leq x_2, \quad 0 \leq x_3 \leq 3
\end{array}
\]
Each constraint in this problem is interpreted by GloptiPoly 3
as a support constraint on the measure associated with variable
$x$, see Section \ref{supcon}:
\begin{verbatim}
>> mpol x 3
>> x(1)+x(2)+x(3) <= 4
Scalar measure support inequality
x(1)+x(2)+x(3) <= 4
\end{verbatim}
The whole problem can be entered as follows:
\begin{verbatim}
>> mpol x 3
>> g0 = -2*x(1)+x(2)-x(3);
>> K = [24-20*x(1)+9*x(2)-13*x(3)+4*x(1)^2-4*x(1)*x(2) ...
 +4*x(1)*x(3)+2*x(2)^2-2*x(2)*x(3)+2*x(3)^2 >= 0, ...
 x(1)+x(2)+x(3) <= 4, 3*x(2)+x(3) <= 6, ...
 0 <= x(1), x(1) <= 2, 0 <= x(2), 0 <= x(3), x(3) <= 3];
>> P = msdp(min(g0), K)
...
Moment SDP problem
  Measure label             = 1
  Relaxation order          = 1
  Decision variables        = 9
  Linear inequalities       = 8
  Semidefinite inequalities = 4x4
\end{verbatim}
The moment problem can then be solved:
\begin{verbatim}
>> [status,obj] = msol(P)
GloptiPoly 3.0
Solve moment SDP problem
...
Global optimality cannot be ensured
status =
     0
obj =
    -6.0000
\end{verbatim}
Since {\tt status=0} the moment SDP problem can be solved
but it is impossible to detect global optimality. The
value {\tt obj=-6.0000} is then a lower bound on the
global minimum of the quadratic problem.

The measure associated with the problem variables
can be retrieved as follows:
\begin{verbatim}
>> mu = meas
Measure 1 on 3 variables: x(1),x(2),x(3)
  with moments of degree up to 2
\end{verbatim}
Its vector of moments can be built as follows:
\begin{verbatim}
>> mv = mvec(mu)
10-by-1 moment vector
(1,1):I[1]d[1]
(2,1):I[x(1)]d[1]
(3,1):I[x(2)]d[1]
(4,1):I[x(3)]d[1]
(5,1):I[x(1)^2]d[1]
(6,1):I[x(1)x(2)]d[1]
(7,1):I[x(1)x(3)]d[1]
(8,1):I[x(2)^2]d[1]
(9,1):I[x(2)x(3)]d[1]
(10,1):I[x(3)^2]d[1]
\end{verbatim}
These moments are the decision variables of the SDP
problem solved with the above {\tt msol} command.
Their numerical values can be retrieved as follows: 
\begin{verbatim}
>> double(mv)
ans =
    1.0000
    2.0000
   -0.0000
    2.0000
    7.6106
    1.4671
    2.3363
    4.8335
    0.5008
    8.7247
\end{verbatim}
The numerical moment matrix can be obtained
using the following commands:
\begin{verbatim}
>> double(mmat(mu))
ans =
    1.0000    2.0000   -0.0000    2.0000
    2.0000    7.6106    1.4671    2.3363
   -0.0000    1.4671    4.8335    0.5008
    2.0000    2.3363    0.5008    8.7247
\end{verbatim}

As explained in \cite{lasserre}, we can build a hierarchy
of nested moment SDP problems, or relaxations, whose
solutions converge monotically and asymptotically to the global optimum,
under mild technical assumptions. By default the command
{\tt msdp} builds the relaxation of lowest order, equal
to half the degree of the highest degree monomial in the 
polynomial data. An additional input argument can be specified
to build higher order relaxations:
\begin{verbatim}
>> P = msdp(min(g0), K, 2)
...
Moment SDP problem
  Measure label             = 1
  Relaxation order          = 2
  Decision variables        = 34
  Semidefinite inequalities = 10x10+8x(4x4)
>> [status,obj] = msol(P)
...
Global optimality cannot be ensured
status =
    0
obj =
   -5.6922
>> P = msdp(min(g0), K, 3)
...
Moment SDP problem
  Measure label             = 1
  Relaxation order          = 3
  Decision variables        = 83
  Semidefinite inequalities = 20x20+8x(10x10)
>> [status,obj] = msol(P)
...
Global optimality cannot be ensured
status =
     0
obj =
    -4.0684
\end{verbatim}
We observe that the moment SDP problems feature an increasing
number of variables and constraints. They generate a mononotically
increasing sequence of lower bounds on the global optimum,
which is eventually reached numerically
at the fourth relaxation:
\begin{verbatim}
>> P = msdp(min(g0), K, 4)
...
Moment SDP problem
  Measure label             = 1
  Relaxation order          = 4
  Decision variables        = 164
  Semidefinite inequalities = 35x35+8x(20x20)
>> [status,obj] = msol(P)
...
2 globally optimal solutions extracted
Global optimality certified numerically
status =
    1
obj =
   -4.0000
>> double(x)
ans(:,:,1) =
    2.0000
    0.0000
    0.0000
ans(:,:,2) =
    0.5000
    0.0000
    3.0000
>> double(g0)
ans(:,:,1) =
   -4.0000
ans(:,:,2) =
   -4.0000
\end{verbatim}

\subsubsection{Rational minimization}

Minimization of a rational function can also be formulated
as a linear moment problem. Given two polynomials $g_0(x)$ and $h_0(x)$,
consider the rational optimization problem
\[
\min_{x \in {\mathbb K}} \frac{g_0(x)}{h_0(x)}
\]
where
\[
{\mathbb K} = \{ x \in {\mathbb R}^n \: :\:
g_i(x) \geq 0, \:i=1,2,\ldots\}
\]
is a basic semialgebraic set described by given
polynomials $g_i(x)$.
Following \cite{jibetean}, the corresponding moment problem is
given by
\[
\begin{array}{ll}
\min_{d\mu \in {\mathbb B}^n} & \int_{\mathbb K} g_0(x)d\mu(x) \\
\mathrm{s.t.} & \int_{\mathbb K} h_0(x)d\mu(x) = 1.
\end{array}
\]
In contrast with the polynomial optimization problem of Section
\ref{constrained}, the optimal measure $d\mu$ supported on $\mathbb K$
is not necessarily a probability measure. Denoting
$h_0(x) = \sum_{\alpha} h_{0\alpha} x^{\alpha}$,
the moments $y_{\alpha}$ of $d\mu$ must satisfy
a linear constraint
\[
\int_{\mathbb K} h_0(x)d\mu(x) = \sum_{\alpha} h_{0\alpha} y_{\alpha} = 1.
\]

As an example, consider the one-variable rational minimization problem
\cite[Ex. 2]{jibetean}:
\[
\min \frac{x^2-x}{x^2+2x+1}.
\]
We can solve this problem with GloptiPoly 3 as follows:
\begin{verbatim}
>> mpol x
>> g0 = x^2-2*x; h0 = x^2+2*x+1;
>> P = msdp(min(g0), mom(h0) == 1);
>> [status,obj] = msol(P)
...
Global optimality certified numerically
status =
     1
obj =
    -0.3333
>> double(x)
ans =
     0.4999
\end{verbatim}

\subsubsection{Several measures}

GloptiPoly 3 can handle several measures whose moments are
linearly related.

For example, consider the GPM arising when solving polynomial
optimal control problems as detailed
in \cite{optcon}. We are seeking two occupation
measures $d\mu_1(x,u)$ and $d\mu_2(x)$ of a 
state vector $x(t)$ and input vector $u(t)$ whose
time variation are governed by the differential equation
\[
\frac{dx(t)}{dt} = f(x,u), \: x(0) = x_0, \: u(0) = u_0
\]
with $f(x,u)$ a given polynomial mapping and $x_0$, $u_0$
given initial conditions. Measure $d\mu_1$ is 
supported on a given semialgebraic set ${\mathbb K}_1$
corresponding to constraints on $x$ and $u$.
Measure $d\mu_2$ is supported on a given semialgebraic set
${\mathbb K}_2$ corresponding to performance requirements.
For example ${\mathbb K}_2 = 0$ indicates that state $x$
must reach the origin.

Given a polynomial test function $g(x)$ we can relax
the dynamics constraint with the moment constraint
\[
\int_{\mathbb{K}_2}g(x)d\mu_2(x)
- g(x_0) = \int_{\mathbb{K}_1}\frac{dg(x)}{dx}f(x,u)d\mu_1(x,u)
\]
linking linearly moments of $d\mu_1$ and $d\mu_2$.
As explained in \cite{optcon}, a lower bound on the minimum
time achievable by any feedback control law $u(x)$ is then
obtained by minimizing the mass of $d\mu_1$ over all possible
measures $d\mu_1$, $d\mu_2$ satisfying the support and moment
constraints. The gap between the lower bound and the exact
minimum time is narrowed by enlarging the class of test
functions $g$.

In the following script we solve this moment
problem in the case of a double integrator
with state and input constraints:
\begin{verbatim}
% bounds on minimal achievable time for optimal control of
% double integrator with state and input constraints

x0 = [1; 1]; u0 = 0; % initial conditions
d = 6; % maximum degree of test function

% analytic minimum time
if x0(1) >= -(x0(2)^2-2)/2
 tmin = 1+x0(1)+x0(2)+x0(2)^2/2;
elseif x0(1) >= -x0(2)^2/2*sign(x0(2))
 tmin = 2*sqrt(x0(1)+x0(2)^2/2)+x0(2);
else
 tmin = 2*sqrt(-x0(1)+x0(2)^2/2)-x0(2);
end

% occupation measure for constraints
mpol x1 2
mpol u1
m1 = meas([x1;u1]);

% occupation measure for performance
mpol x2 2
m2 = meas(x2);

% dynamics
scaling = tmin; % time scaling
f = scaling*[x1(2);u1];

% test function
g1 = mmon(x1,d);
g2 = mmon(x2,d);

% initial condition
assign([x1;u1],[x0;u0]);
g0 = double(g1);

% moment problem
P = msdp(min(mass(m1)),...
 u1^2 <= 1,... % input constraint
 x1(2) >= -1,... % state constraint
 x2'*x2 <= 0,... % performance = reach the origin
 mom(g2) - g0 == mom(diff(g1,x1)*f)); % linear moment constraints

% solve
[status,obj] = msol(P);
obj = scaling*obj;

disp(['Minimum time = ' num2str(tmin)]);
disp(['LMI ' int2str(d) ' lower bound = ' num2str(obj)])
\end{verbatim}
For the initial condition $x_0=[1 \:\: 1]$ the exact
minimum time is equal to $3.5$. In Table \ref{bounds} we
report the monotically increasing sequence of lower bounds
obtained by solving moment problems with test functions
of increasing degrees. We used the above script and the
semidefinite solver SeDuMi 1.1R3.

\begin{table}[h]
\begin{center}
\begin{tabular}{c|cccccccc}
degree & 2 & 4 & 6 & 8 & 10 & 12 & 14 & 16 \\ \hline
bound & 1.0019 & 2.3700 & 2.5640 & 2.9941 & 3.3635 & 3.4813 & 3.4964 & 3.4991
\end{tabular}
\caption{Minimum time optimal control for double integrator with state
and input constraints: lower bounds on exact minimal time $3.5$ achieved by
solving moment problems with test functions of increasing degrees.
\label{bounds}}
\end{center}
\end{table}

\subsection{Using YALMIP}\label{yalmip}

By default GloptiPoly 3 uses the semidefinite solver SeDuMi \cite{sedumi}
for solving numerically SDP moment problems. It is however possible
to use any solver interfaced through YALMIP \cite{yalmip}
by setting a configuration flag with the {\tt mset} command:
\begin{verbatim}
>> mset('yalmip',true)
\end{verbatim}
Parameters for YALMIP, handled with the YALMIP command {\tt sdpsettings},
can be forwarded to GloptiPoly 3 with the {\tt mset} command. For example,
the following command tells YALMIP
to use the SDPT3 solver (instead of SeDuMi)
when solving moment problems with GloptiPoly:
\begin{verbatim}
>> mset(sdpsettings('solver','sdpt3'));
\end{verbatim}

\subsection{SeDuMi parameters settings}\label{settings}

The default parameters settings of SeDuMi \cite{sedumi}
can be altered as follows:
\begin{verbatim}
>> pars.eps = 1e-10;
>> mset(pars)
\end{verbatim}
where {\tt pars} is a structure of parameters
consistent with SeDuMi's format.

\subsection{Exporting moment SDP problems}

A moment problem {\tt P} of class {\tt msdp}
can be converted into SeDuMi's input format:
\begin{verbatim}
>> [A,b,c,K] = msedumi(P);
\end{verbatim}
The SDP problem can then be solved with SeDuMi as follows:
\begin{verbatim}
>> [x,y,info] = sedumi(A,b,c,K);
\end{verbatim}
See \cite{sedumi} for more information on SeDuMi's input data format.

Similarly, a moment SDP problem can be converted into YALMIP's input format:
\begin{verbatim}
>> [F,h,y] = myalmip(P);
\end{verbatim}
where variable {\tt F} contains the LMI constraints (YALMIP class {\tt lmi}), {\tt h}
is the objective function (YALMIP class {\tt sdpvar}) and {\tt y} is
the vector of moments (YALMIP class {\tt sdpvar}).
The SDP problem can then be solved with any semidefinite
solver interfaced through
YALMIP as follows:
\begin{verbatim}
>> solvesdp(F,h);
>> ysol = double(y);
\end{verbatim}

\subsection{Moment substitutions}

By performing explicit moment substitutions it is often possible
to reduce significantly
the number of variables and constraints in moment SDP problems.
Version 2 of GloptiPoly implemented these substitutions for mixed-integer
0-1 problems only \cite{gloptipoly2}. With version 3, these substitutions
can be carried out in full generality.

GloptiPoly 3 carries out moment substitutions as soon as the left hand-side
of a support or moment equality constraint consists of an isolated monic monomial.
Otherwise, no substitution is achieved and the equality constraint is
preserved.

For example, consider the $AW_2^9$ Max-Cut problem studied in \cite[\S 4.7]{gloptipoly2},
with variables $x_i$ taking values $-1$ or $+1$ for $i=1,\ldots,9$. These integer
constraints can be expressed algebraically as $x^2_i=1$. The following piece of
code builds up the third relaxation of this problem:
\begin{verbatim}
>> W = diag(ones(8,1),1)+diag(ones(7,1),2)+diag([1 1],7)+diag(1,8);
>> W = W+W'; n = size(W,1); e = ones(1,n); Q = (diag(e*W)-W)/4;
>> mset clear
>> mpol('x', n)
>> P = msdp(max(x'*Q*x), x.^2 == 1, 3)
GloptiPoly 3.0
Define moment SDP problem
  Valid objective function
  Number of support constraints = 9 including 9 substitutions
  Number of moment constraints = 0
Measure #1
  Maximum degree = 2
  Number of variables = 9
  Number of moments = 5005
Order of SDP relaxation = 3
Mass of measure 1 set to one
Total number of monomials = 5005
Perform moment substitutions
Perform support substitutions
Number of monomials after substitution = 465
Generate moment and support constraints
Generate moment SDP problem  

Moment SDP problem 
  Measure label             = 1
  Relaxation order          = 3
  Decision variables        = 465
  Semidefinite inequalities = 130x130
\end{verbatim}
We see that out of the 5005 moments (corresponding to all the monomials
of 9 variables of degree up to 6), only 465 linearly independent moments
appear in a reduced moment matrix of dimension 130.

With the following syntax, moment substitutions are not carried out:
\begin{verbatim}
>> P = msdp(max(x'*Q*x), x.^2-1 == 0, 3)
...
Mass of measure 1 set to one
Total number of monomials = 5005
Perform moment substitutions
Number of monomials after substitution = 5004
Generate moment and support constraints
Generate moment SDP problem 

Moment SDP problem 
  Measure label             = 1
  Relaxation order          = 3
  Decision variables        = 5004
  Linear equalities         = 6435
  Semidefinite inequalities = 220x220
\end{verbatim}
Only the mass is substituted, and the remaining 5004 moments
linked by 6435 linear equalities (many of which are redundant)
now appear explicitly in a full-size moment matrix of dimension 220.


\begin{thebibliography}{10}

\bibitem{akhiezer}
N. I.  Akhiezer. The classical moment problem. Hafner, New York, 1965.

\bibitem{akhiezer2}
N. I. Akhiezer, M. G. Krein. Some questions in the theory of
moments. American Mathematical Society Translations Vol. 2, 1962.

\bibitem{gloptipoly2}
D. Henrion, J. B. Lasserre. GloptiPoly: global optimization over
polynomials with Matlab and SeDuMi. ACM Transactions on Mathematical
Software, Vol. 29, No. 2, pp. 165-194, 2003.

\bibitem{jibetean}
D. Jibetean, E. de Klerk. Global optimization of rational functions: a
semidefinite programming approach. Mathematical Programming, Vol. 106,
No. 1, pp. 93-109, 2006.

\bibitem{landau}
H. J. Landau. Moments in mathematics. In H. J. Landau (Editor), Proceedings
of Symposium on Applied Mathematics, Vol. 37, American Mathematical Society, 1980.

\bibitem{lasserre}
J. B. Lasserre. Global optimization with polynomials and the problem
of moments. SIAM Journal on Optimization, Vol. 11, No. 3,
pp. 796--817, 2001.

\bibitem{optcon}
J. B. Lasserre, C. Prieur, D. Henrion. Nonlinear optimal control: numerical
approximation via moments and LMI relaxations. Proceedings of the joint
IEEE Conference on Decision and Control and European Control Conference,
Sevilla, Spain, December 2005.

\bibitem{gpm}
J. B. Lasserre. A semidefinite programming approach to the generalized
problem of moments. To appear in Mathematical Programming, Series B, 2007.

\bibitem{laurent}
M. Laurent. Moment matrices and optimization over polynomials - A survey
on selected topics. CWI Amsterdam, The Netherlands, September 2005. 

\bibitem{yalmip}
J. L\"ofberg. YALMIP : a toolbox for modeling and optimization in Matlab.
Proceedings of the IEEE Symposium on Computer-Aided Control System Design (CACSD),
Taipei, Taiwan, 2004. See {\tt control.ee.ethz.ch/$\sim$joloef/yalmip.php}

\bibitem{sedumi}
J. F. Sturm and the Advanced Optimization Laboratory at McMaster
University, Canada. SeDuMi version 1.1R3, October 2006. See
{\tt sedumi.mcmaster.ca}

\end{thebibliography}
\end{document}